\documentclass[a4paper]{article}
\usepackage{amsmath}
\newcommand{\bea}{\begin{eqnarray*}}
\newcommand{\eea}{\end{eqnarray*}}

\newcommand{\D}{\cal D}
\newcommand{\T}{{\cal T}}

\newcommand{\la}{\langle}
\newcommand{\ra}{\rangle}
\newcommand{\Rm}{I\!\!R^m}

\begin{document}

\title{Equivariant Operators between some\\ Modules of the Lie Algebra of Vector Fields\footnote{This work was supported by
MCESR Grant MEN/CUL/99/007. The author thanks P. Lecomte and P.
Mathonet for helpful comments.}}
\author{N. Poncin}\maketitle

\newtheorem{re}{Remark}
\newtheorem{rem}{Provisional remark}
\newtheorem{theo}{Theorem}
\newtheorem{prop}{Proposition}
\newtheorem{lem}{Lemma}
\newtheorem{cor}{Corollary}

\begin{abstract}
The space ${\cal D}_p^k$ of differential operators of order $\le
k$, from the differential forms of degree $p$ of a smooth manifold
$M$ into the functions of $M,$ is a module over the Lie algebra of
vector fields of $M,$ when it's equipped with the natural Lie
derivative. In this paper, we compute all equivariant i.e.
intertwining operators $T:{\cal D}_p^k \rightarrow {\cal
D}_q^\ell$ and conclude that the preceding modules of differential
operators are never isomorphic. We also answer a question of P.
Lecomte, who observed that the restriction of some homotopy
operator---introduced in \cite{PLc}---to ${\cal D}_p^k$ is
equivariant for small values of $k$ and $p$.\end{abstract}



\section{Introduction}

Let $M$ be a smooth, Hausdorff, second countable, connected
manifold of dimension $m$.

Denote by $\Omega^p(M)$ the space of differential forms of degree
$p$ of $M$, by $N$ the space $C^\infty(M)$ of smooth functions of
$M,$ and by ${\cal D}_p^k$ the space of differential operators of
order smaller than or equal to $k,$ from $\Omega^p(M)$ into $N$.
If $X \in Vect(M),$ where $Vect(M)$ is the Lie algebra of vector
fields of $M,$ and $D \in {\cal D}_p^k$, the Lie derivative
\begin{equation} L_X D = L_X  \circ D - D \circ L_X \label{lie}\end{equation}
is a differential operator of order at most $k$ and so $({\cal
D}_p^k,L)$ is a $Vect(M)$-module.

In this paper, we shall determine all the spaces ${\cal
I}_{p,q}^{k,\ell}$ of \textit{equivariant} operators from ${\cal
D}_p^k$ into ${\cal D}_q^\ell$, that is all operators $T:{\cal
D}_p^k \rightarrow {\cal D}_q^\ell$ such that
\begin{equation}
L_X \circ T = T \circ L_X, \forall X \in Vect(M) \label{equi}.
\end{equation}
In \cite{PLc2}, P. Lecomte, P. Mathonet, and E. Tousset computed
all linear equivariant mappings between modules of differential
operators acting on densities. We solve here an analogous problem,
thus answering a question of P. Lecomte, who noticed that some
homotopy operator, which locally coincides---up to a
coefficient---with the Koszul differential (see \cite{PLc}), is
equivariant if it's restricted to low order differential operators
and asked whether this property still holds for higher orders.

\section{Local representation}\label{locrepr}

The formalism described in this section was also used in
\cite{NP1} and \cite{NP2}.

Consider an open subset $U$ of $I\!\!R^m$, two real
finite-dimensional vector spaces $E$ and $F$, and some local
operator
\[O\in {\cal L}(C^\infty(U,E),C^\infty(U,F))_{loc}.\]
The operator is fully defined by its values on the products $fe$,
$f\in C^{\infty}(U), e\in E$. A well known theorem of J. Peetre
(see \cite{JP}) states that it has the form
\[O(fe)=\sum_{\alpha}O_{\alpha}(\partial^{\alpha}(fe))
=\sum_{\alpha}O_{\alpha}(e)\partial^{\alpha}f,\] where
$\partial_x^{\alpha}=\partial_{x^1}^{\alpha^1}\ldots\partial_{x^m}^{\alpha^m}$
and $O_{\alpha}\in C^{\infty}(U,{\cal L}(E,F))$. Moreover, the
coefficients $O_{\alpha}$ are well determined by $O$ and the
series is locally finite (it is finite, if $U$ is relatively
compact).

We shall symbolize the partial derivative $\partial^{\alpha}f$ by
the monomial
$\xi^{\alpha}=\xi_1^{\alpha^1}\ldots\linebreak\xi_m^{\alpha^m}$ in
the components $\xi_1,\ldots,\xi_m$ of some linear form $\xi\in
(I\!\!R^m)^*$, or---at least mentally---even by $\xi^{\alpha}f,$
if this is necessary to avoid confusion. The operator $O$ is thus
represented by the polynomial
\[{\cal O}(\xi;e)=\sum_{\alpha}O_{\alpha}(e)\xi^{\alpha}.\]
When identifying the space $Pol (I\!\!R^m)^*$ of polynomials on
$(I\!\!R^m)^*$ with the space $\vee I\!\!R^m$ of symmetric
contravariant tensors of $I\!\!R^m$, one has ${\cal O}\in
C^{\infty}(U,\vee I\!\!R^m\otimes {\cal L}(E,F))$. Let's emphasize
that the form $\xi$ symbolizes the derivatives in $O$ that act on
the argument $fe\in C^{\infty}(U,E)$, while $e\in E$ represents
this argument. In the sequel, we shall no longer use different
notations for the operator $O$ and its representative polynomial
${\cal O}$; in order to simplify notations, it's helpful to use
even the same typographical sign, when referring to the argument
$fe$ and its representation $e$.

Let's for instance look for the local representation of
(\ref{lie}). If $D\in {\cal D}_p^k$, its restriction $D\!\!\mid_U$
(or simply $D,$ if no confusion is possible) to a domain $U$ of
local coordinates of $M$, is a local operator from $C^\infty (U,
\wedge ^p (I\!\!R^m )^*)$ into $C^\infty(U)$ that is represented
by $D(\omega)\simeq D(\xi;\omega)$, where $\omega\in C^\infty (U,
\wedge ^p (I\!\!R^m )^*)$ on the l.h.s. and $\omega\in\wedge ^p
(I\!\!R^m )^*$ on the r.h.s., and where $\xi$ represents the
derivatives acting on $\omega$. The Lie derivative of $D(\omega)$
with respect to a vector field $X\in C^\infty (U;I\!\!R^m)$, is
then represented by $L_X(D(\omega))\simeq \la X,\eta+\xi\ra
D(\xi;\omega)$. Here $\eta\in (I\!\!R^m)^*$ is associated to $D$
and $\la X,\eta+\xi\ra$ denotes the evaluation of $X\in I\!\!R^m$
on $\eta +\xi$. When associating $\zeta$ to $X$, one gets
$D(L_X\omega)\simeq D(\xi+\zeta;\la X,\xi\ra\omega +\zeta\wedge
i_X\omega)$ and
\begin{equation}
(L_X D)(\omega) \simeq \la X,\eta\ra D(\xi;\omega) - \la X,\xi\ra
\tau_\zeta D(\xi;\omega) - D(\xi + \zeta; \zeta \wedge i_X
\omega),\label{Lie}\end{equation} where $\tau_\zeta
D(\xi;\omega)=D(\xi+\zeta;\omega)-D(\xi;\omega)$.

\section{Locality of equivariant operators}

Let $U$ be a domain of local coordinates of $M$. If $D \in {\cal
D}_p^k\!\!\mid_U$ (${\cal D}_p^k\!\!\mid_U$ is defined similarly
to ${\cal D}_p^k,$ but for $M=U$), its representation is a
polynomial $D$ in $C^\infty (U,E_p^k)$, where $E_p^k = \vee ^{\le
k} I\!\!R^m \otimes \wedge ^p I\!\!R^m$, or---equivalently---in
the space of smooth sections $\Gamma(\vee ^{\le k}TU\otimes\wedge
^pTU)$. We denote by $\sigma(D)$ or simply $\sigma$ the
homogeneous part of maximal degree of this polynomial i.e. the
principal symbol of the considered differential operator. If $L^t$
is the natural Lie derivative of tensor fields, one easily
verifies that
\[\left(L_X^t(\sigma(D))\right)(\xi;\omega)=\la X,\eta\ra\sigma(\xi;\omega)
-\la X,\xi\ra
(\zeta\partial_{\xi})\sigma(\xi;\omega)-\sigma(\xi;\zeta\wedge
i_X\omega),\] where $\eta$ and $\zeta$ are associated to $\sigma$
and $X$ respectively, and where $\zeta\partial_{\xi}$ denotes the
derivation with respect to $\xi$ in the direction of $\zeta$. If
$L^{op}$ is the formerly defined Lie derivative $L$ of
differential operators, it follows from (\ref{Lie}) that
$\left(\sigma(L^{op}_XD)\right)(\xi;\omega)=\left(L_X^t(\sigma(D))\right)(\xi;\omega)$
i.e. that the principal symbol is equivariant with respect to all
vector fields.\\

We are now prepared to prove the following lemma:

\begin{lem} Every equivariant operator $T\in {\cal
I}_{p,q}^{k,\ell}$ is local.\end{lem}

\textit{Proof.} It suffices to show that the family ${\cal
L}_p^k=\{L_X:{\cal D}_p^k\longrightarrow {\cal D}_p^k\mid X\in
Vect(M)\}$ is closed with respect to locally finite sums and is
locally transitive (lt); this means that each point of $M$ has
some neighborhood $\Omega,$ such that for every open subset
$\omega\subset\Omega$ and every $D\in {\cal D}_p^k $ with compact
support in $\omega$, $D$ can be decomposed into Lie derivatives,
\begin{equation}D=\sum_{i=1}^n L_{X_i}D_i \;\;\; (D_i\in {\cal D}_p^k,\mbox{
supp }\!X_i\subset\omega,\mbox{ supp
}\!D_i\subset\omega)\label{loctrans}, \end{equation} with $n$ independent of
$D$ and $\Omega$. Indeed, proposition $3$ of \cite{DWLc} states
that ${\cal L}_p^k$ is then globally transitive i.e. that
(\ref{loctrans}) holds for every open subset $\omega$ of $M$ and
every $D\in {\cal D}_p^k$ with support in $\omega$. Take now a
differential operator $D\in {\cal D}_p^k$ that vanishes in an open
subset $V$ of $M$ and let $x_0$ be an arbitrary point of $V$.
Since supp $\!D\subset\omega=M\backslash\omega_0$ for some
neighborhood $\omega_0$ of $x_0$, $D$ has form (\ref{loctrans})
and $T(D)(x_0)=\sum_{i=1}^n \left(L_{X_i}(T(D_i))\right)(x_0)=0$.

In order to confirm local transitivity (LT), let's point out that
example 12 of \cite{DWLc} shows that the family
${\pounds}_s^r=\{L_X:\Gamma(\otimes_s^rM)\longrightarrow\Gamma(\otimes_s^rM)\mid
X\in Vect(M)\}$ is lt, if \begin{equation}s-r\neq m. \label{condloc}\end{equation}
Moreover, the domains $U$ of the charts of $M$ that correspond to
cubes in $I\!\!R^m$, play the role of the neighborhoods $\Omega$
in the definition of LT.

LT of ${\cal L}_p^k$ may now be proved as follows. If $x\in M$,
take $\Omega =U$ and if $D\in {\cal D}_p^k$ is a differential
operator with compact support in some open subset $\omega$ of $U$,
note that its principal symbol $\sigma(D)$ is a contravariant
tensor field with compact support in $\omega$. Hence condition
(\ref{condloc}) is satisfied and
\[\sigma(D)=\sum_{i=1}^nL_{X_i}^t\sigma(D_i)=\sigma\left(\sum_{i=1}^nL_{X_i}^{op}D_i\right)\]
($D_i\in {\cal D}_p^k$; supp $\!X_i$, supp $\!D_i\subset\omega$;
$n$ independent of $D$ and $U$). Thus
$D-\sum_{i=1}^nL_{X_i}^{op}D_i\in {\cal D}_p^{k-1}$ and we
conclude by iteration.  \rule{1.5mm}{2.5mm}

\section{Local expression of the equivariance equation}

Let $U$ be a connected, relatively compact domain of local
coordinates of $M$.

Recall that if $D \in {\cal D}_p^k\!\!\mid_U$, its representation
is a polynomial $D \in C^\infty (U,E_p^k)$, where $E_p^k = \vee
^{\le k} I\!\!R^m \otimes \wedge ^p I\!\!R^m$.

We identify ${\cal D}_p^k\!\!\mid_U$ with $C^\infty (U, E_p^k)$.
Thus, if $T \in {\cal I}_{p,q}^{k,\ell}$, the restriction
$T\!\!\mid _U$ is a local operator from $C^\infty (U, E_p^k)$ into
$C^\infty (U, E_q^\ell)$, with representation
$T(\eta;D)(\xi;\omega)$ $(\eta, \xi \in (I\!\!R^m)^*, D \in E_p^k,
\omega \in \wedge ^q (I\!\!R^m )^*)$.

It's easily checked that equivariance condition (\ref{equi})
locally reads
\begin{eqnarray}
\lefteqn{(X.T)(\eta ;D)(\xi ;\omega ) -  \la X,\eta\ra (\tau
_\zeta
T(\eta ;D))(\xi ;\omega )}\nonumber\\
 & & -  \la X,\xi\ra \tau _\zeta  (T(\eta;D))(\xi ;\omega) + T(\eta  + \zeta ;X\tau _\zeta  D)(\xi
 ;\omega)\nonumber\\
 & & - T(\eta;D)(\xi  + \zeta ;\zeta \wedge i_X \omega ) + T(\eta  + \zeta ;D( \cdot  + \zeta ;\zeta  \wedge i_X  \cdot))(\xi;\omega)=0,
\label{fonda}\end{eqnarray} where $X.T$ is obtained by derivation
of the coefficients of $T$ in the direction of $X$.

Take in equation (\ref{fonda}) the terms of degree 0 in $\zeta$:
\[(X.T)(\eta;D)(\xi;\omega) = 0.\]
This means that the coefficients of $T$ are constant.

The terms of degree $1$ lead to the equation
\begin{eqnarray*}
\lefteqn{\la X,\eta\ra (\zeta \partial_\eta) T(\eta;D)(\xi;\omega)}\\
 & & - T(\eta;X(\zeta\partial_\xi)D)(\xi;\omega)-T(\eta;D(\cdot\,;\zeta \wedge i_X \cdot))(\xi;\omega)\\
 & & +\la X,\xi\ra(\zeta\partial_\xi)T(\eta;D)(\xi;\omega)+ T(\eta;D)(\xi;\zeta \wedge i_X\omega)= 0,
\end{eqnarray*} which, if $\rho$ denotes the natural action of $gl(m,I\!\!R)$, may be rewritten
\[\rho(X \otimes \zeta) \left( T(\eta;D)(\xi;\omega) \right) = 0.\]
Note that $T(\eta;D)(\xi;\omega)$ is completely characterized by
$T(\eta; Y^r \otimes (X_1 \wedge \ldots \wedge X_p))(\xi; \nu^1
\wedge \ldots \wedge \nu^q)$ $(Y,X_i \in \Rm, \, \nu^j \in
(\Rm)^*, \, r \in \{0, \ldots, k\})$. This last expression is a
polynomial in $Y,X_i$ and in $\eta,\xi,\nu^j$. It thus follows
from the description of invariant polynomials under the action of
$gl(m,I\!\!R)$ (see \cite{HW}), that it is a polynomial ${\cal
T}_r(\la Y,\xi\ra,\la Y,\eta\ra,\la Y,\nu^j\ra,\la X_i,\xi\ra,\la
X_i,\eta\ra,\la X_i,\nu^j\ra)$ in the evaluations of the vectors
on the linear forms.

In order to determine the most general structure of ${\cal T}_r$,
observe that this polynomial is homogeneous of degree $r$ in $Y$
and degree $1$ in the $X_i$'s and the $\nu^j$'s, and that
furthermore it's antisymmetric in the $X_i$'s and the $\nu^j$'s.
It follows from the skew-symmetry in the $\nu^j$'s, that $Y$ is
evaluated on at most one $\nu^j$, so that $q\le p+1,$ and from the
skew-symmetry in the $X_i$'s, that $\xi$ and $\eta$ are evaluated
on at most one $X_i$, so that $q\ge p-2$. Finally, $q$ is $p-2,
\,p-1, \,p$ or $p+1$. Set now $\Lambda = X_1 \wedge \ldots \wedge
X_p$, $\omega = \nu^1 \wedge \ldots \wedge \nu^q$, $u = \la
Y,\xi\ra$ and $v = \la Y,\eta\ra$. The following possible forms of
the terms of ${\cal T}_r$ and the corresponding conditions on $p$
and $q$ are immediate consequences of the preceding observations.
\begin{equation}
\begin{array}{llll}
    \text{term type}&&\text{condition}&\text{term form}\\\\
    \text{(1) no } \la Y,\nu^j\ra&&\\\\
    \text{(1.1) no } \la X_i,\xi\ra&\text{ no } \la X_i,\eta\ra&q=p&v^su^{r-s}\la\Lambda,\omega\ra\\
    \text{(1.2) one } \la X_i,\xi\ra&\text{ no }
    \la X_i,\eta\ra&q=p-1&v^su^{r-s}\la i_\xi\Lambda,\omega\ra\\
    \text{(1.3) no } \la X_i,\xi\ra&\text{ one }
    \la X_i,\eta\ra&q=p-1&v^su^{r-s}\la i_\eta\Lambda,\omega\ra\\
    \text{(1.4) one } \la X_i,\xi\ra&\text{ one }
    \la X_i,\eta\ra&q=p-2&v^su^{r-s}\la i_{\xi}i_\eta\Lambda,\omega\ra\\\\
    \text{(2) one }  \la Y,\nu^j\ra&&\\\\
    \text{(2.1) no } \la X_i,\xi\ra&\text{ no } \la X_i,\eta\ra&q=p+1&v^su^{r-s-1}\la Y\wedge\Lambda,\omega\ra\\
    \text{(2.2) one } \la X_i,\xi\ra&\text{ no }
    \la X_i,\eta\ra&q=p&v^su^{r-s-1}\la Y \wedge i_\xi\Lambda,\omega\ra\\
    \text{(2.3) no } \la X_i,\xi\ra&\text{ one }
    \la X_i,\eta\ra&q=p&v^su^{r-s-1}\la Y \wedge i_\eta\Lambda,\omega\ra\\
    \text{(2.4) one } \la X_i,\xi\ra&\text{ one }
    \la X_i,\eta\ra&q=p-1&v^su^{r-s-1}\la Y\wedge i_{\xi}i_\eta\Lambda,\omega\ra\\
\end{array}
\label{form}
\end{equation}
\begin{re} These structures of ${\cal T}_r$ show that a priori
equivariant operators are mappings $T:{\cal D}_p^k \longrightarrow
{\cal D}_{p-2}^{k+1},\,T:{\cal D}_p^k \longrightarrow {\cal
D}_{p-1}^{k+1},\,T:{\cal D}_p^k \longrightarrow {\cal D}_p^k\mbox{
or } T:{\cal D}_p^k \longrightarrow {\cal
D}_{p+1}^{k-1}.$\label{refoaprio}\end{re}

It's interesting to rewrite (\ref{fonda}) in terms of the ${\cal
T}_r$'s $(r \in \{0, \ldots, k\})$. To simplify notations, we set
\[\begin{array}{llll}
    \lambda = \la X,\xi\ra &\mu = \la X,\eta\ra &\nu = \la X,\zeta\ra &\pi^j =
    \la X,\nu^j\ra\\
    u = \alpha_0 = \la Y,\xi\ra &v = \beta_0 = \la Y,\eta\ra &\gamma_0 = \la Y,\zeta\ra &\delta_0^j = \la
    Y,\nu^j\ra\\
    \alpha_i = \la X_i,\xi\ra &\beta_i = \la X_i,\eta\ra &\gamma_i = \la X_i,\zeta\ra &\delta_i^j =
    \la X_i,\nu^j\ra
\end{array}\]
Substitute in (\ref{fonda}), $Y^r \otimes (X_1 \wedge \ldots
\wedge X_p)$ and $\nu^1 \wedge \ldots \wedge \nu^q$ to $D$ and
$\omega$ respectively, and use subscript $\iota\in\{0, \ldots,
p\}$. Note that:

\begin{itemize}
\item
Taylor expansion gives
\[T(\eta;D)(\xi +\zeta;\omega)=\sum\frac{1}{a!}(\zeta\partial_{\xi})^{a}T(\eta;D)(\xi;\omega)
=\sum\frac{1}{a!}(\zeta\partial_{\xi})^{a}{\cal
T}_r(\alpha_\iota,\beta_\iota,\delta_\iota^j)\]
\item
the computation of the successive directional derivatives
$(\zeta\partial_{\xi})^{a}$ of the composite function ${\cal
T}_r(\alpha_\iota,\beta_\iota,\delta_\iota^j)$ leads to the same
rule than the computation of the successive powers of a sum of
$p+1$ real terms:
\begin{eqnarray*}
\lefteqn{(\zeta\partial_{\xi})^{a}{\cal
T}_r(\alpha_\iota,\beta_\iota,\delta_\iota^j)}\\
 & & = \sum_{\varrho_0 + \ldots + \varrho_p = a} \; \frac{a!}{\varrho_0 ! \ldots \varrho_p !} \;
 \gamma_0^{\varrho_0} \ldots \gamma_p^{\varrho_p}\left(\partial_{\alpha_0}^{\varrho_0} \ldots
 \partial_{\alpha_p}^{\varrho_p}{\cal T}_r\right)
 (\alpha_\iota,\beta_\iota,\delta_\iota^j)
\end{eqnarray*}
\item
\[\tau_{\zeta}Y^r=\sum_{a=0}^{r-1}\left(\begin{array}{c}r\\a\end{array}\right)\gamma_0^{r-a}Y^{a}\]
\item
\[XY^{a}=\frac{1}{a+1}(X\partial_Y)Y^{a+1}\]
\item
\begin{eqnarray} \lefteqn{(X_1\ \wedge \ldots \wedge
X_p)(\zeta \wedge i_X \cdot)}\nonumber\\
 & & = X \wedge i_\zeta(X_1\ \wedge
\ldots \wedge X_p)\nonumber\\
 & & = \sum_{b=1}^p \gamma_b (X
\partial_{X_b})(X_1\ \wedge \ldots \wedge X_p)
\label{hint}\end{eqnarray}
\end{itemize}
Transform now the terms of (\ref{fonda}) in conformity with the
preceding hints. If $\left(\partial{\cal
T}\right)(\cdot)=\left(\partial{\cal
T}\right)(\alpha_{\iota},\cdot,\delta_{\iota}^j)$ and if
$d_r(\xi)$ denotes the degree of ${\cal T}_r$ in $\xi$, the local
equivariance equation (\ref{fonda}) finally reads in terms of the
${\cal T}_{r}$'s:
\begin{eqnarray}
\lefteqn{\mbox{\boldmath $\mu$} \left({\cal T}_r(\beta_{\iota} +
\mbox{\boldmath $\gamma_\iota$}) - {\cal T}_r(\beta_\iota)
\right)}\nonumber\\
 & & +\mbox{\boldmath $\lambda$} \sum_{a=1}^{d_r(\xi)} \;
\; \sum_{\varrho_0 + \ldots + \varrho_p = a} \; \frac{1}{\varrho_0
! \ldots \varrho_p !} \; \mbox{\boldmath $\gamma_0^{\varrho_0}
\ldots \gamma_p^{\varrho_p}$}
\left(\partial_{\alpha_0}^{\varrho_0} \ldots
\partial_{\alpha_p}^{\varrho_p}{\cal
T}_r\right)(\beta_\iota)\nonumber\\
 & & + \sum_{b=1}^q \; \mbox{\boldmath $\pi^b$}
\sum_{a=0}^{d_{r}(\xi)} \; \sum_{\iota = 0}^p \mbox{\boldmath
$\gamma_\iota$}\nonumber\\
 & & \sum_{\varrho_0 + \ldots + \varrho_p = a} \;
\frac{1}{\varrho_0 ! \ldots \varrho_p !} \; \mbox{\boldmath
$\gamma_0^{\varrho_0} \ldots \gamma_p^{\varrho_p}$}
\left(\partial_{\alpha_0}^{\varrho_0} \ldots
\partial_{\alpha_p}^{\varrho_p}
\partial_{\delta_\iota^b}{\cal
T}_r\right)(\beta_\iota)\label{Tr}\\
 & & - \sum_{a=0}^{r-1} \; \left(\begin{array}{c}r\\a
\end{array}\right) \frac{1}{a+1} \; \mbox{\boldmath
$\gamma_0^{r-a}$}
\begin{array}[t]{l}
\lefteqn{\left[ \right. \mbox{\boldmath
$\lambda$}\left(\partial_{\alpha_0} {\cal
T}_{a+1}\right)(\beta_\iota + \mbox{\boldmath
$\gamma_\iota$})}\\
+ \mbox{\boldmath $(\mu + \nu)$}\left(\partial_{\beta_0} {\cal
T}_{a+1}\right)(\beta_\iota + \mbox{\boldmath
$\gamma_\iota$})\\
+ \; \sum_{j=1}^q \; \mbox{\boldmath
$\pi^j$}\left(\partial_{\delta_0^j} {\cal
T}_{a+1}\right)(\beta_\iota + \mbox{\boldmath
$\gamma_\iota$})\left.\right]
\end{array}\nonumber\\
 & & - \sum_{b=1}^p \; \mbox{\boldmath $\gamma_b$}
\sum_{a=0}^r \; \left(\begin{array}{c} r\\a \end{array} \right)
\mbox{\boldmath $\gamma_0^{r-a}$} \;
\begin{array}[t]{l}
\lefteqn{\left[ \right. \mbox{\boldmath
$\lambda$}\left(\partial_{\alpha_b} {\cal
T}_{a}\right)(\beta_\iota + \mbox{\boldmath
$\gamma_\iota$})}\\
+ \mbox{\boldmath $(\mu + \nu)$}\left(\partial_{\beta_b} {\cal
T}_{a}\right)(\beta_\iota + \mbox{\boldmath
$\gamma_\iota$})\\
+ \; \sum_{j=1}^q \; \mbox{\boldmath
$\pi^j$}\left(\partial_{\delta_b^j} {\cal
T}_{a}\right)(\beta_\iota + \mbox{\boldmath
$\gamma_\iota$})\left.\right] = 0,
\end{array}\nonumber
\end{eqnarray}
for each $r \in \{0, \ldots, k\}$.\\

\begin{re}
In the sequel, we suppose that dimension $m$ is not only $\ge
sup(p,q)$, but also $\ge inf\{p+2, q+3\}$. In the skipped
borderline cases, results have to be completed, but proofs are not
elegant.
\end{re}

We claim that (\ref{Tr}) is a polynomial identity in the
independent variables $\lambda, \ldots, \delta_p^q$. Indeed, since
as well the vectors $X,Y,X_i$ as the forms $\xi,\eta,\zeta,\nu^j$
are arbitrary, take $X=e_1,Y=e_2,X_i = e_{i+2}$ ($e_k$: canonical
basis of $\Rm$) and let the $p+2$ first components of the
preceding forms in the dual basis of $e_k$ vary in $I\!\!R$, if
$inf\{p+2, q+3\}=p+2$; proceed similarly but
exchange roles of vectors and forms, if $inf\{p+2, q+3\}=q+3$.\\

When seeking in (\ref{Tr}) the terms of degree $1$ in $\nu$ and
$\gamma_0$, and of degree $0$ in $\gamma_1, \ldots, \gamma_p$ (in
the sequel we shall denote these terms by $(\nu)^1 \gamma_0^1
\gamma_1^0 \ldots \gamma_p^0$), one gets
\begin{equation} \partial_{\beta_0} {\cal T}_r = 0, \label{beta0}
\end{equation}
where ${\cal T}_0$ is a priori independent of $Y$ and $\beta_0 =
\la Y,\eta\ra$. The terms in $(\nu)^1 \gamma_0^0
\gamma_1^0\linebreak\ldots \gamma_i^1 \ldots \gamma_p^0$ ($i \in
\{1, \ldots, p\}, \: p
> 0$), $(\lambda)^1 \gamma_0^2 \gamma_1^0 \ldots \gamma_p^0$, $(\lambda)^1 \gamma_0^1 \gamma_1^0 \ldots \gamma_i^1
\ldots \gamma_p^0$ ($i \in \{1, \ldots, p\}, \: p > 0$),
$(\pi^j)^1 \gamma_0^2 \gamma_1^0 \ldots \gamma_p^0$ $(j \in \{1,
\ldots, q\}, \: q > 0)$ and $(\pi^j)^1 \gamma_0^1 \gamma_1^0
\ldots \gamma_i^1\linebreak\ldots\gamma_p^0$ \, $(i \in \{1,
\ldots, p\}, \: j \in \{1, \ldots, q\}, \: p > 0, \: q > 0)$,
read:
\begin{equation}
\partial_{\beta_i} {\cal T}_r = 0, \label{betai}
\end{equation}
\begin{equation}
\partial_{\alpha_0}^2 {\cal T}_r - r
\partial_{\alpha_0} {\cal T}_{r-1} = 0, \label{lambda00}
\end{equation}
\begin{equation}
\partial_{\alpha_0 \alpha_i}{\cal T}_r - r \partial_{\alpha_i} {\cal T}_{r-1}=
0,
\label{lambda0i}
\end{equation}
\begin{equation}
2 \partial_{\alpha_0 \delta_0^j}{\cal T}_r - r
\partial_{\delta_0^j}{\cal T}_{r-1} = 0
\label{pi00}
\end{equation}
resp.
\begin{equation}
\partial_{\alpha_i \delta_0^j} {\cal T}_r + \partial_{\alpha_0
\delta_i^j} {\cal T}_r - r
\partial_{\delta_i^j} {\cal T}_{r-1} = 0.
\label{pi0i}
\end{equation}
These partial equations will allow to compute all equivariant
operators.

\section{Determination of the equivariant operators}

\begin{prop}
Equivariant operators $T\in {\cal I}_{p,q}^{k,\ell}$ are mappings
\[
T:{\cal D}_p^k \longrightarrow {\cal D}_{p-1}^{k+1}, \; T:{\cal
D}_p^k \longrightarrow {\cal D}_p^k \mbox{  or  } T:{\cal D}_p^k
\longrightarrow {\cal D}_{p+1}^{k-1}.
\]
\end{prop}

\textit{Proof.} If $T\in {\cal I}_{p,p-2}^{k,k+1}$, it follows
from (\ref{form}) and (\ref{beta0}) that ${\cal T}_r = c_r
\alpha_0^r \det(\alpha_i, \beta_i, \delta_i^j)$ ($r\in
\{0,\ldots,k\},c_r\in I\!\!R$) and from (\ref{betai}) that ${\cal
T}_r = 0$ ($r\in \{0,\ldots,k\}$). Hence the result (see remark
\ref{refoaprio}). \rule{1.5mm}{2.5mm}

\subsection{Case \boldmath{$q = p+1$}}

\begin{prop}
All spaces ${\cal I}_{p,p+1}^{k,k-1}$ vanish, except\\\\
(i) the spaces ${\cal I}_{p,p+1}^{1,0}$ with bases defined by equation (\ref{inv11}), and\\
(ii) the space ${\cal I}_{0,1}^{2,1}$ with basis defined by
equation (\ref{inv12}).
\end{prop}

\textit{Proof.} Equations (\ref{form}) and (\ref{beta0}) show that
${\cal T}_r = c_r \alpha_0^{r-1}\det(\delta_{\iota}^j)\:(r \in
\{0, \ldots, k\},\linebreak c_r\in I\!\!R).$ Moreover, equations
(\ref{lambda00}) and (\ref{pi00}) yield \begin{equation}(r-1)c_r = rc_{r-1} \;
(r \in \{3, \ldots, k\}),\label{c1}\end{equation}
\begin{equation}
2(r-1)c_r = rc_{r-1} \; (r \in \{2, \ldots, k\}).
\label{c2}\end{equation} If $k\ge 3$, (\ref{c1}) and (\ref{c2})
imply that each invariant vanishes. If $k=2$, (\ref{c2}) confirms
that $c_1 = c_2 = c$ ($c\in I\!\!R$). If in addition $p>0$,
(\ref{pi0i}), written for $r=2,$ gives $c=0$. If $(k,p)=(2,0)$ or
$k=1$, all local invariants have the form $\T_0 = 0, \T_1 = c\la
Y, \nu \ra, \T_2 = c\la Y, \xi \ra \la Y, \nu \ra$ resp. $\T_0 =
0, \T_1 = c\det(\delta_{\iota}^j)$.\\

In \cite{PLc}, P. Lecomte introduced---in a more general
framework---some homotopy operator $K$ for the dual $d^*$ of the
de Rham differential $d$.

It's well-known that each $D'\in{\D}^1_p$ admits a global
decomposition $D' = \sum \langle \Lambda, L_X \cdot \rangle + \sum
\langle \Omega, \cdot \rangle$, where the sums are finite, where
$\Lambda$ and $\Omega$ are antisymmetric contravariant
$p\,$-tensors on $M,$ and where $X$ denotes a vector field of $M.$
Similarly, any differential operator $D''\in{\D}^2_0$ may be
written $D'' = \sum \Lambda \: L_X \circ L_Y  +  \sum \Omega \:
L_Z +\sum\Theta,$ with $\Lambda,\Omega,\Theta\in N$ and $X,Y,Z \in
Vect(M)$.

One easily verifies that \begin{equation}K\!\!\mid_{{\D}^1_p}(D') =
\frac{1}{1+p}\sum \langle \Lambda, i_X \cdot \rangle \label{inv11}
\end{equation} and that \begin{equation}K\!\!\mid_{{\D}^2_0}(D'') = \frac{1}{2} \sum
\Lambda \; (i_XL_Y + L_Xi_Y) + \sum \Omega \; i_Z.\label{inv12}\end{equation}
For the first formula for instance, it suffices to remember that
\[ \langle \Lambda, L_X\omega \rangle \simeq \langle X,\xi \rangle
(X_1 \wedge \ldots \wedge X_p)(\omega)+ (X_1 \wedge \ldots \wedge
X_p)(\zeta \wedge i_X\omega),\] with unmistakable notations (see
section \ref{locrepr}), and to use the local representation of $K$
(see \cite{PLc}).

These results imply that---in the implicated cases---the homotopy
operator is independent of the local coordinates and the partition
of unity involved in its construction. Furthermore, the r.h.s. of
(\ref{inv11}) and (\ref{inv12}) is independent of the chosen
decomposition of $D'$ and $D''$ respectively. It's now obvious
that
\[K:{\D}^1_p\rightarrow {\D}^0_{p+1}\mbox{   and   }
K:{\D}^2_0\rightarrow {\D}^1_1\] are equivariant operators.

Finally, the spaces ${\cal I}^{1,0}_{p,p+1}$ and ${\cal
I}^{2,1}_{0,1}$ are generated by $K$. \rule{1.5mm}{2.5mm}

\subsection{Case \boldmath{$q = p-1$  $(p>0)$}}

\begin{prop}
The spaces ${\cal I}_{p,p-1}^{k,k+1}$ are one-dimensional vector
spaces with basis $d^*$. \label{th-1}\end{prop}

\textit{Proof.} Equations (\ref{form}), (\ref{beta0}), and
(\ref{betai}) tell us that $\T_r = c_r \alpha_0^r \det(
\alpha_i,\delta_i^j)$ $(r\in\{0, \ldots, k \},c_r\in I\!\!R)$ and
equation (\ref{lambda0i}) brings out that $c_r = c$ $(r\in\{0,
\ldots, k\},c\in I\!\!R)$.\\

Since obviously $d^*\in {\cal I}^{k,k+1}_{p,p-1}$, the conclusion
follows. \rule{1.5mm}{2.5mm}

\subsection{Case \boldmath{$q = p$}}

\begin{prop}
The spaces ${\cal I}_{p,p}^{k,k}$ are one-dimensional, except
${\cal I}_{0,0}^{k,k}$ $(k>0)$ and ${\cal I}_{p,p}^{1,1}$ $(p>0)$
that have dimension $2$. Possible bases are $id$, $(id,I_0)$ resp.
$(id, d^*K)$, where $I_0$ is defined by $I_0: D \rightarrow
D(1)id$ $(1$ stands for the constant function $x \rightarrow 1)$.
\end{prop}

\textit{Proof.} (i) Look first at the case $p=0$. Equations
(\ref{form}), (\ref{beta0}), and (\ref{lambda00}) show that
$\T_0=c_0\;(c_0\in I\!\!R)$ and $\T_r=c\alpha_0^r\;(r \in
\{1,\ldots,k\},c\in I\!\!R).$\\

If $k=0$, the space of invariants is generated by $id$. Otherwise,
dimension of ${\cal I}^{k,k}_{0,0}$ is $2$ and the invariants
$id$ and $I_0$ form a basis.\\

(ii) If $p>0$, it follows from equations (\ref{form}),
(\ref{hint}), (\ref{beta0}), and (\ref{betai}) that \[ \T_r =
c_r\alpha_0^r\det(\delta_i^j)+d_r\alpha_0^{r-1}\sum_{n=1}^p \:
\alpha_n\mbox{$\det_{n0}$}(\delta_i^j)\;\;(r \in
\{0,\ldots,k\},c_r,d_r\in I\!\!R).\] In this expression,
$\det_{n0}(\delta_i^j)$ denotes the determinant
$\det(\delta_i^j)$, where the line
$(\delta_n^1,\ldots,\delta_n^p)$ has been replaced by
$(\delta_0^1,\ldots,\delta_0^p)$. When exploiting equations
(\ref{lambda0i}) and (\ref{pi00}), you find
$(r-1)d_r=rd_{r-1}\;\;(r \in \{2,\ldots,k\})$ resp.
$2(r-1)d_r=rd_{r-1}\;\;(r \in \{2,\ldots,k\})$, so that $d_r=0$
($r\in\{0,\ldots,k\}$), if $k\ge 2$. Apply now equation
(\ref{pi0i}). If $k\ge 2$, this condition shows that \[ \T_r =
c\alpha_0^r\det(\delta_i^j)\;\;(r \in \{0,\ldots,k\},c\in
I\!\!R)\] and if $k=1$, it entails that
$c_1-c_0+d_1=0$.\\

If $k=0$ or $k\ge 2$, invariants are thus generated by $id$. In
the case $k=1$, dimension is $2$ and $id,d^*K\in {\cal
I}^{1,1}_{p,p}$ are possible generators. \rule{1.5mm}{2.5mm}
\begin{re}
Equation (\ref{inv11}) reveals that
$(d^*K)(D)=(1/(1+p))\sum\la\Lambda,$ $i_X d\cdot\ra$, if
$D=\sum\la\Lambda,L_X\cdot\ra+\sum\la\Omega,\cdot\ra$.\end{re}

\noindent Universit\'e de Luxembourg\\D\'epartement de Math\'ematiques\\Avenue de la Fa\"{\i}encerie, 162 A\\
L-1511 Luxembourg, Luxembourg\\E-mail: poncin@cu.lu
\end{document}